# Decomposing a Multi-Scale Optimization Framework for Grid-Integrated Electrolysis using Aggregate-Informed Benders

**Kiernan X. Jennings[a], Victor M. Zavala[ab], and Styliana Avraamidou[a]***

[a] University of Wisconsin-Madison, Department of Chemical and Biological Engineering, Madison, WI, USA
[b] Mathematics and Computer Science Division, Argonne National Laboratory, Lemont, IL, USA
* Corresponding Author: avraamidou@wisc.edu

## ABSTRACT

Demand response (DR) operation of electrolysis devices is gaining traction to capitalize on volatile electricity markets, but their dynamic operation poses challenges to the durability and lifespan of these systems. Our recently developed multi-scale optimization framework for grid-integrated electrolysis systems studied the impacts of DR and effects on device durability. A major hurdle in that work is tractably scaling the model to include participation in high-frequency electricity markets and/or longer time horizons. To address this, in this work, we developed an aggregate-informed Benders decomposition method to tractably solve large instances of this problem structure. To illustrate the efficacy, we solve a case study with market participation in the day-ahead market (DAM) and real-time markets (RTM) for up to 40 years to effectively capture the effects of device lifespan change. We compare our decomposition algorithm to both commercial solvers and traditional Benders decomposition. We find that incorporating aggregate subproblem information accelerates convergence, reducing the final optimality gap by up to ~81% relative to traditional Benders on the 40-year horizon instances that otherwise stall near 85%. We additionally apply the algorithm to find that RTM participation offers economic advantages under the higher price volatility.

## INTRODUCTION

The introduction of volatile and intermittent energy supply and demand technologies causes issues for the electricity grid; rapid fluctuation of energy availability requires more flexible technologies to counteract these swings. Wholesale electricity markets hedge against this uncertainty by adopting a hierarchical market structure with different markets, including the day-ahead market (DAM) and the real-time market (RTM). This structure enables power consumers to participate under a demand response (DR) regime to provide flexibility and perform arbitrage on electricity prices [1].

Electrolysis devices offer a synergistic opportunity to accelerate the decarbonization of the power grid and chemical manufacturing. Through DR, these devices can provide flexibility and reduce operational expenses relative to constant operation schemes [2]. One drawback of this is that the rapid ramping of power in these devices has been shown in short- to medium-term studies to adversely affect the durability of these systems [3]. A mixed-integer linear programming (MILP) framework developed by Jennings et al. [4] modeled the multi-scale coupling of optimal DR operation with degradation and electrolysis stack replacement scheduling. This framework enabled studying market participation over a 22-year time horizon while participating in the DAM. However, tractably solving high-frequency market participation at similar time horizons or longer remained a major computational challenge.

Decomposition algorithms are a well-known technique that can address these tractability issues. The general concept is to divide large, high-dimensional optimization problems into lower-dimensional subproblems and coordinate information between them. One widely known exact decomposition algorithm is Benders decomposition [5], which exploits the block-diagonal problem structure present in the modeling framework of [4]. This algorithm partitions the model into a master problem and subproblem(s), where complicating integer variables are fixed so the subproblem can be solved, and its solution generates cutting planes to refine the master problem. However, a drawback with traditional Benders decomposition is that the master problem lacks subproblem information until cuts accumulate, so early solutions are poor and thus convergence is limited.

In this work, we present an aggregate-informed Benders decomposition algorithm for solving the multi-scale MILP electrolysis system market participation framework. We generate a valid aggregate-informed master problem to produce stronger first-stage candidate solutions. Then, we rigorously prove the validity of this approach. Furthermore, we demonstrate the application of this decomposition algorithm through a case study on the DR operation scheduling of an industrial alkaline water electrolysis (AWE) device across multiple time horizons and two market types (DAM & RTM).

## OPTIMIZATION FRAMEWORK

Jennings et al. [4] presented a multi-scale MILP optimization framework of electrolysis systems for replacement and operation scheduling when participating in electricity markets. This monolithic problem structure can be mathematically represented as:

$$\begin{aligned} \min_{z,x} \quad & c^\top z + \sum_{m \in M} d_m^\top x_m \\ \text{s.t.} \quad & w^\top z \le r \\ & A_m x_m + G_m z_m \le b_m, \quad \forall m \in \mathcal{M} \\ & L_{m-1} x_{m-1} - L_m x_m = 0, \quad \forall m \in \mathcal{M} \setminus \{m_1\} \\ & z \in \{0,1\}^N \\ & x_m \in \{0,1\}^{4|\mathcal{D}\times\mathcal{J}|} \times \mathbb{R}^{4|\mathcal{D}\times\mathcal{J}|} \quad \forall m \in \mathcal{M} \end{aligned} \tag{1}$$

Here, we are minimizing the negative net present value (NPV) of the electrolysis system over the given time horizon of $N$ years. The stack replacement decisions, $z$, affect the DR operation scheduling decisions, $x$, for each year $m$. The linking constraint represents the degrading efficiency state of the electrolysis system over time. The reader is referred to [4] for additional details on the formulation.

## DECOMPOSITION FRAMEWORK

We use the market participation model in (1) to hierarchically decompose the replacement and operation problems.

### Benders Decomposition

In Benders decomposition, the problem (1) is partitioned into two: the master problem (2) and the subproblem (3). Here, we use the multi-cut variation of Benders where we assign a separate variable to each subproblem [6].

$$\begin{aligned} \min_{z,\theta} \quad & c^\top z + \sum_{m \in \mathcal{M}} \theta_m \\ \text{s.t} \quad & w^\top z \le r \\ & \theta_m \ge \{\text{cuts}\} \quad \forall m \in \mathcal{M} \\ & z \in \{0,1\}^N \end{aligned} \tag{2}$$

Where $z_m$ represents stack replacement decisions in year $m$ and the cost-to-go $\theta_m$ is the approximation of the yearly subproblem. The cost-to-go variable underestimates the optimal value of the subproblem and is iteratively refined from the generation of cutting planes from the optimal subproblem solution. We adapt multiple cutting planes given the binary decision variables in the subproblem. We use a combination of Benders optimality cuts (from the LP relaxation), L-shaped optimality cuts [7] and simple integer feasibility cuts.

The subproblem can be represented as:

$$\begin{aligned} \Omega_m(z_m, \bar{L}_{m-1}) = \min_{x} \quad & d_m^\top x_m \\ \text{s.t} \quad & A_m x_m \le b_m - G_m \bar{z}_m \\ & L_m x_m = \bar{L}_{m-1} \\ & x_m \in \{0,1\}^{4|\mathcal{D}\times\mathcal{J}|} \times \mathbb{R}^{4|\mathcal{D}\times\mathcal{J}|} \end{aligned} \tag{3}$$

where the binary hourly operation decisions, $x$, are constrained by hydrogen production, electricity consumption, daily hydrogen demands, and operation logic.

For the subproblem solve, we fix the master problem's stack replacement decision variables $z$ and solve the annual subproblem (3) for a parameterized initial efficiency state of the electrolyzer, $\bar{L}_{m-1}$. We solve each yearly subproblem sequentially to recover the upper bound (UB) of the program. This greedy approach reproduces a feasible solution for the operational subproblems.

### Aggregate Benders Decomposition

We propose to aggregate the annual subproblems to produce a tractable under-approximation of the cost-to-go variable in the master problem. This produces a lifted-form instance of the full-space monolithic problem whilst maintaining the constraint structure [8]. This can be adapted to introduce valid supporting cutting planes for the cost-to-go variable, essentially warm-starting Benders [9].

To perform the aggregation procedure, we partition the time horizon of a subproblem for a given year $m$ into the chronological set of clusters $C_m$. These clusters of a given year define aggregate variables of the full-space problem, which is a valid relaxation of the feasible space [8].

$$\sum_{t \in c} x_t = \hat{x}_c \qquad \forall c \in C_m$$

This partitioning procedure for the $m$ subproblem is abstract to any scheme which maintains the chronological ordering of the clusters. With this, we generate the ag-gregate-informed master problem:

$$\begin{aligned}\min_{z,\hat{x},\theta} \quad & c^\top z + \sum_{m\in\mathcal{M}} \theta_m \\ \text{s.t} \quad & w^\top z \le r \\ & \hat{A}_c \hat{x}_c + G_c z \le \hat{b}_c \quad \forall c \in \mathcal{C}_m \\ & \hat{L}_{c-1}\hat{x}_{c-1} + \hat{L}_c \hat{x}_c = 0 \quad \forall c \in \mathcal{C}_m \setminus \{c_1\} \\ & \theta_m \ge \sum_{c\in\mathcal{C}_m} \hat{d}_c^\top \hat{x}_c \quad \forall m \\ & \theta_m \ge \{\text{cuts}\}_m \\ & z \in \{0,1\}^N, \hat{x}_c \in \mathbb{R}^{N|\mathcal{C}_m|}\end{aligned} \tag{4}$$

where $\theta_m$ is under-approximated by the aggregate subproblem objective function.

## Valid Lower Bound

The aggregate variables defined previously must produce a valid relaxation on the original problem to rigorously bound the cost-to-go variable. It remains nontrivial to produce such relaxation of the full-space parameters given the time varying electricity market price data.

As such, we seek an under-approximation of the electricity market prices for a given cluster $c$ to relax the operation scheduling subproblem in year $m$. The aggregation procedure is performed *a priori* by partitioning the historical electricity market price data. We define the minimum price point as the constant price for the given cluster:

$$\hat{\psi}_c^E = \min_{t\in c} \psi_t^E \quad \forall c \in \mathcal{C}_m$$

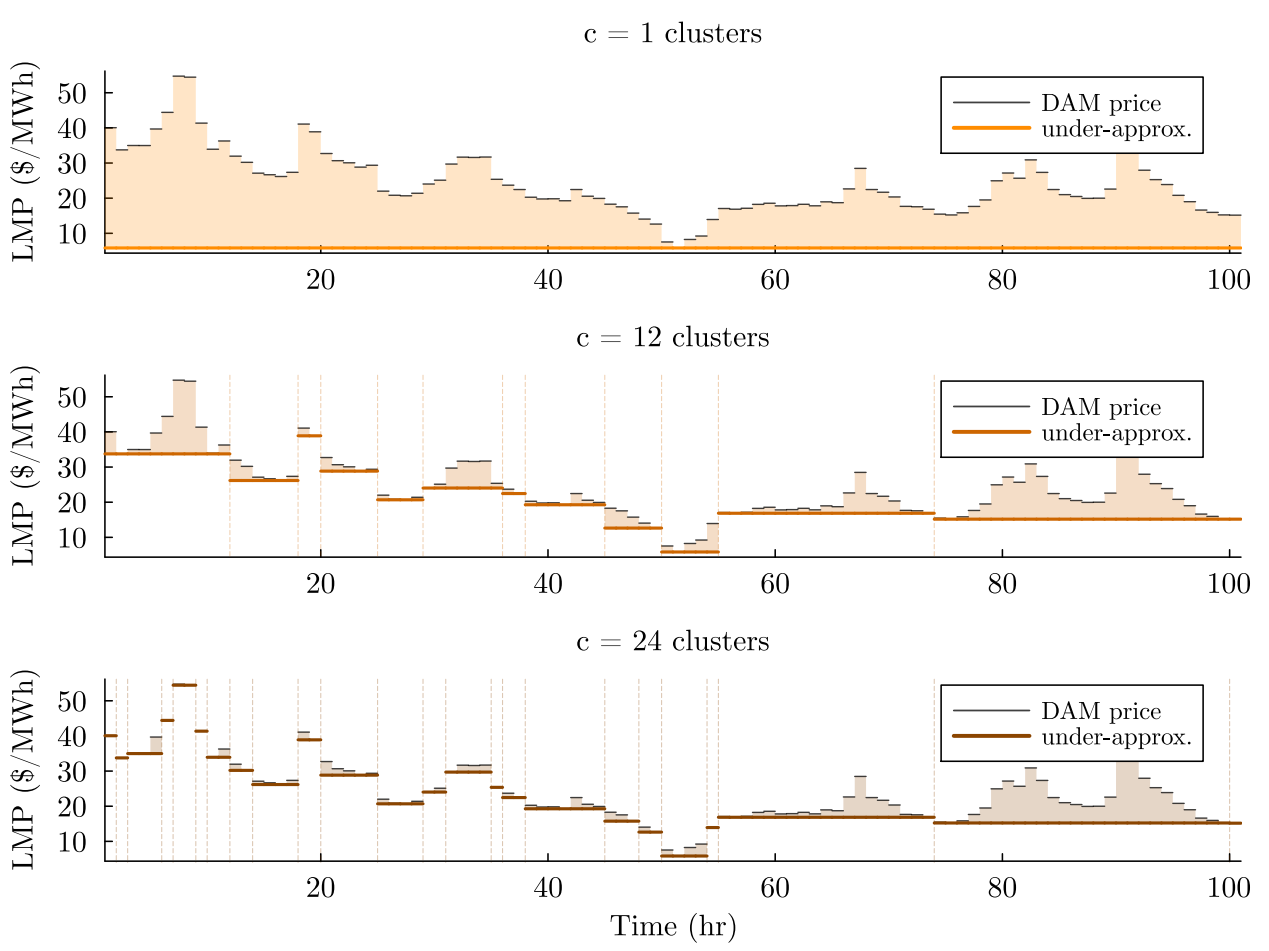

**Figure 1.** DAM electricity price aggregation for various number of clusters on representative 100-hour sample.

We use a heuristic algorithm to produce these clusters which aims to minimize the error between the under approximator and the true market price. Figure 1 displays the varying levels of clustering for a representative 100-hour DAM price series. As we define higher levels of aggregation, $c$, the problem approaches the full monolithic problem ($c = 8760$ for 1 year of hourly DAM prices).

From this definition of $\hat{\psi}_c^E$, we prove the validity of this lower bounding procedure and apply this algorithm in practice.

**Theorem 1.**

Let $\Omega_m = d_m^\top x_m$ be the optimal objective value of subproblem $m$, and let

$$\hat{\Omega}_m = \sum_{c\in\mathcal{C}_m} \hat{d}_c^\top \hat{x}_c$$

be the optimal objective value of the aggregate subproblem $m$ over a set of clusters $\mathcal{C}_m$. If the aggregate parameter $\hat{\psi}_c^E$ is defined as the local minimum price within each cluster, then the aggregate subproblem produces a valid lower bound on the original Benders subproblem, such that $\hat{\Omega}_m \le \Omega_m$.

**Proof 1.**

By definition, $\hat{\psi}_c^E \le \psi_t^E \;\; \forall t \in c$ where $c \in \mathcal{C}_m$. For any nonnegative power consumption $e_t$ and summing over the hierarchical time sets,

$$\hat{\psi}_c^E \sum_{t\in c} e_t \le \sum_{t\in c} \psi_t^E e_t \quad \forall c \in \mathcal{C}_m$$

Subtracting the hydrogen revenue terms and using the aggregate terms for electricity expenses and hydrogen revenue, respectively, $\hat{e}_c = \sum_{t\in c} e_t$, $\hat{h}_c = \sum_{t\in c} h_t$

$$\hat{\psi}_c^E \hat{e}_c - \psi^H \hat{h}_c \le \sum_{t\in c} (\psi_t^E e_t - \psi^H h_t) \; \forall c \in \mathcal{C}_m$$

Both sides are multiplied by the non-negative annual discount rate,

$$\frac{1}{(1+\rho)^m}\left(\hat{\psi}_c^E \hat{e}_c - \psi^H \hat{h}_c\right) \le \frac{1}{(1+\rho)^m} \sum_{t\in c} (\psi_t^E e_t - \psi^H h_t)$$

Thus, $\hat{\Omega}_m \le \Omega_m$. Q.E.D.

**Algorithm 1** Aggregate-Informed Benders Decomposition

Initialize: $LB \leftarrow -\infty, UB \leftarrow +\infty, \epsilon = 0.01, k \leftarrow 0$
**While** $\frac{|UB-LB|}{UB} \ge \epsilon$ and $k \le k^{max}$
  **Master problem:** Solve aggregate MP (4)
  $\phi = \min_{z,\hat{x},\theta} c^\top z + \sum_{m\in M} \theta_m$
  $LB^{(k+1)} \leftarrow \max(LB^{(k)}, \phi)$
  $\bar{z}^{(k)} \leftarrow \arg\min c^\top z + \sum_{m\in M} \theta_m$
  **Subproblem:** Solve sequential subproblems (3)
  $\Gamma \leftarrow 0$
  **for $m = 1$ to $N$ do**
    $\Omega_m\left(\bar{z}_m^{(k)}, \bar{L}_{m-1}\right) = \min_{x_m} d_m^\top x_m$
    **if** subproblem $m$ is feasible **then**
      $\{cuts\}_m \leftarrow$ optimality cut
    **else**
      $\{cuts\} \leftarrow$ feasibility cut
    **end if**
      $\Gamma \leftarrow \Gamma + \Omega_m$
  **end for**

if all subproblems are feasible then
$UB^{(k+1)} \leftarrow \min\{UB^{(k)}, \Gamma\}$
end if
$k \leftarrow k + 1$
end while

## CASE STUDY: MARKET PARTICIPATION OF AWE

The proposed algorithm is applied to a case study and to assess an electrolyzer device during long-term market participation. Jennings et al. [4] analyzed the DR operation scheduling when participating in the DAM for a 2.2 MW AWE using historical data for Texas electricity markets (ERCOT). We extend this case study by comparing commercial solvers and the proposed decomposition approach in the DAM and the 15-minute RTM for the same electrolyzer specifications.

We use progressively longer time horizons ($N$ = 10, 20, 40 years) and varying levels of aggregation ($|C_m|$ = 1, 12, 24). The monolithic problem sizes are noted in Table 1. To reach the designated time horizon, historical price data is extracted from 2014-2024 and replicated to reach desired time horizon. All results from the presented work are solved using JuMP modeling language [9] and using Intel(R) Xeon(R) CPU E5-2698 v3 @ 2.30GHz with 32 processor cores. All optimization problems were solved using the Gurobi 13.0 solver. The threshold MIP gap on the problem was set to 1%. We set a maximum Benders iteration count of 60. No time limit on the algorithm was set, but a 300s time limit was set on a single annual subproblem. All results can be reproduced using the code provided in https://github.com/Avraamidou-Research-Group-CESE/aggregate-informed-benders.

**Table 1:** Computational statistics for problem instances.

| Statistic | Time Horizon (yr) N = 10 | 20 | 40 |
|---|---|---|---|
| *Day-Ahead Market* | | | |
| Continuous var. | 350,410 | 700,820 | 1,401,640 |
| Binary var. | 350,410 | 700,820 | 1,401,640 |
| Equality constr. | 350,401 | 700,801 | 1,401,601 |
| Inequality constr. | 967,295 | 1,934,595 | 3,869,195 |
| *Real-Time Market* | | | |
| Continuous var. | 1,401,610 | 2,803,220 | 5,606,440 |
| Binary var. | 1,401,610 | 2,803,220 | 5,606,440 |
| Equality constr. | 1,401,600 | 2,803,200 | 5,606,400 |
| Inequality constr. | 3,858,098 | 7,716,198 | 15,432,398 |

Table 2 shows the overall computational results from our aggregate-informed Benders decomposition as well as the traditional (multi-cut) Benders and monolithic problem. For both the DAM and RTM participation scenarios, 5 of 6 monolithic problems were unable to be computed due to out-of-memory (OOM) errors. The one tractable instance, DAM for 10 years, took around 9.2x the amount of time as the respective shortest problem instance.

**Table 2:** Optimality gap (%) and wall-clock time (s) for each algorithm across planning horizons. Bold values indicate the best result within each metric group. The termination criteria of the algorithms are the maximum number of Benders iterations ($k^{max}$=60). *OOM* = out-of-memory termination.

| Problem Instance | | Monolithic | | Traditional Benders | | Aggregate-Informed Benders | | | | | |
|---|---|---|---|---|---|---|---|---|---|---|---|
| | | | | | | $c = 1$ | | $c = 12$ | | $c = 24$ | |
| Electricity Market | Time Horizon (yr) | Gap (%) | Time (s) | Gap (%) | Time (s) | Gap (%) | Time (s) | Gap (%) | Time (s) | Gap (%) | Time (s) |
| Day-Ahead | 10 | 1.00 | 4180 | **0.00** | 1203 | 0.88 | 469 | 0.80 | 466 | 0.76 | **456** |
| | 20 | *OOM* | | 9.15 | 8510 | 1.82 | 7595 | 1.40 | **5611** | **1.30** | 5904 |
| | 40 | *OOM* | | 86.53 | **8301** | 5.64 | 10629 | 4.35 | 10886 | **4.18** | 10886 |
| Real-Time | 10 | *OOM* | | 0.73 | 4240 | 0.33 | 2483 | **0.23** | 1870 | 0.91 | **1470** |
| | 20 | *OOM* | | 11.47 | **7521** | 1.38 | 11850 | 1.21 | 10901 | **1.03** | 10479 |
| | 40 | *OOM* | | 85.44 | **1248** | 17.76 | 9529 | 5.08 | 20194 | **3.02** | 18053 |

As previously mentioned, traditional Benders can produce poor first stage decisions leading to slow convergence. This is especially prevalent in the 40-year instances where the optimality gap remained around 85% for both the DAM and RTM. In both scenarios, traditional Benders produces many infeasible replacement schedules, which induces the generation of weak cutting planes. Subsequently, we reach the iteration limit set on the algorithm with generally poor optimality gap convergence in a short time period.

From Table 2, we see that our algorithm converges faster than traditional Benders. The 40-year DAM participation trial showcases an 80.89% relative optimality gap difference between traditional Benders and $c = 1$ aggregate-informed Benders. At these long time horizons, traditional Benders decomposition is unable to generate feasible replacement schedules to test for subproblems, thus encountering many infeasible subproblem solutions and weak feasibility cut generation. However, the aggregate-informed master problem provides superior first stage variable candidates which, in turn, produces stronger cutting planes from the subproblem in fewer iterations. We see that the inclusion of the aggregate subproblem greatly accelerates algorithm convergence by proposing better candidate decisions in the master problem.

We also see that the increasing level of aggregation coincides with this superior algorithm convergence; the $c = 24$ aggregate-informed Benders outperformed with 4 of 6 instances having the smallest optimality gap. As we approach larger clustering $c$, we increase the information in the master problem thus improving the under-approximation on the cost-to-go.

We also note that the time required to reach termination varied tremendously (ranges up to 18,946s) between the algorithms. We see that the largest problem (40-year RTM) took traditional Benders algorithm a minimum of 1,248s, then 9,529s in the $c = 1$ case and up to 20,194s in the $c = 12$ case for our algorithm. Traditional Benders requires the least time to reach the iteration limit, again due to the poor first stage decisions and weak infeasibility cut generation. For our algorithm, one would expect the higher clustering levels to proportionally increase the computational time required. This is not the case in any of the trials. Future work must look to study the tractability of the master problem as a function of the clustering level.

Figure 2 displays the convergence of the upper and lower bounds for RTM participation over 20 years for the four Benders algorithms. First, we see that the instantiation of the traditional Benders is much faster than our presented aggregate-informed algorithm with the first upper and lower bound generated in 513.8s compared to 1524.1s. We see that the computational burden shifts to more time per iteration but reduces the necessary iterations to converge in our algorithm. Notably, we see that traditional Benders gets to a similar optimality gap when our algorithm gets the first bounds but fails to close the gap in similar time frames. Our aggregate-informed Benders decomposition algorithm demonstrates superior convergence to the optimal solution than the traditional Benders in practice. This is further supported by the final convergence to within 1.38% and as low as 1.03% in Table 2 for 20-year RTM participation.

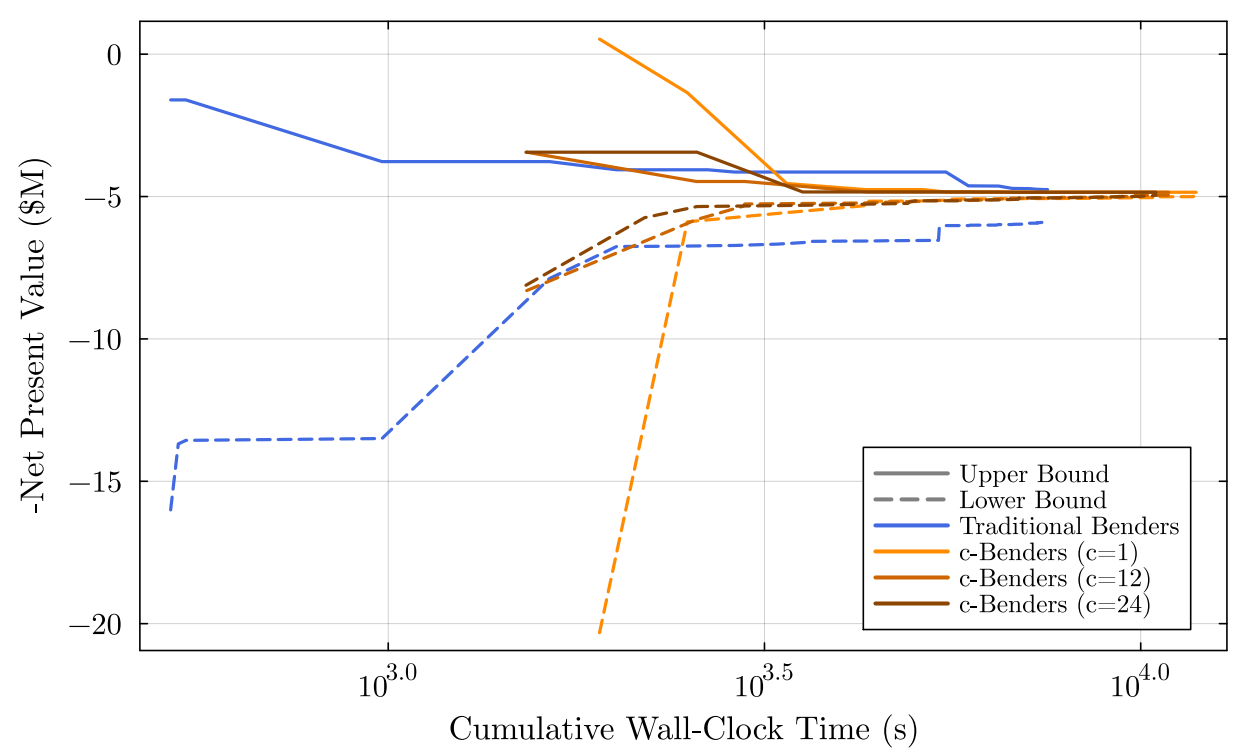

**Figure 2.** Convergence of traditional and aggregate-informed Benders (c-Benders) for RTM participation over a 20-year time horizon.

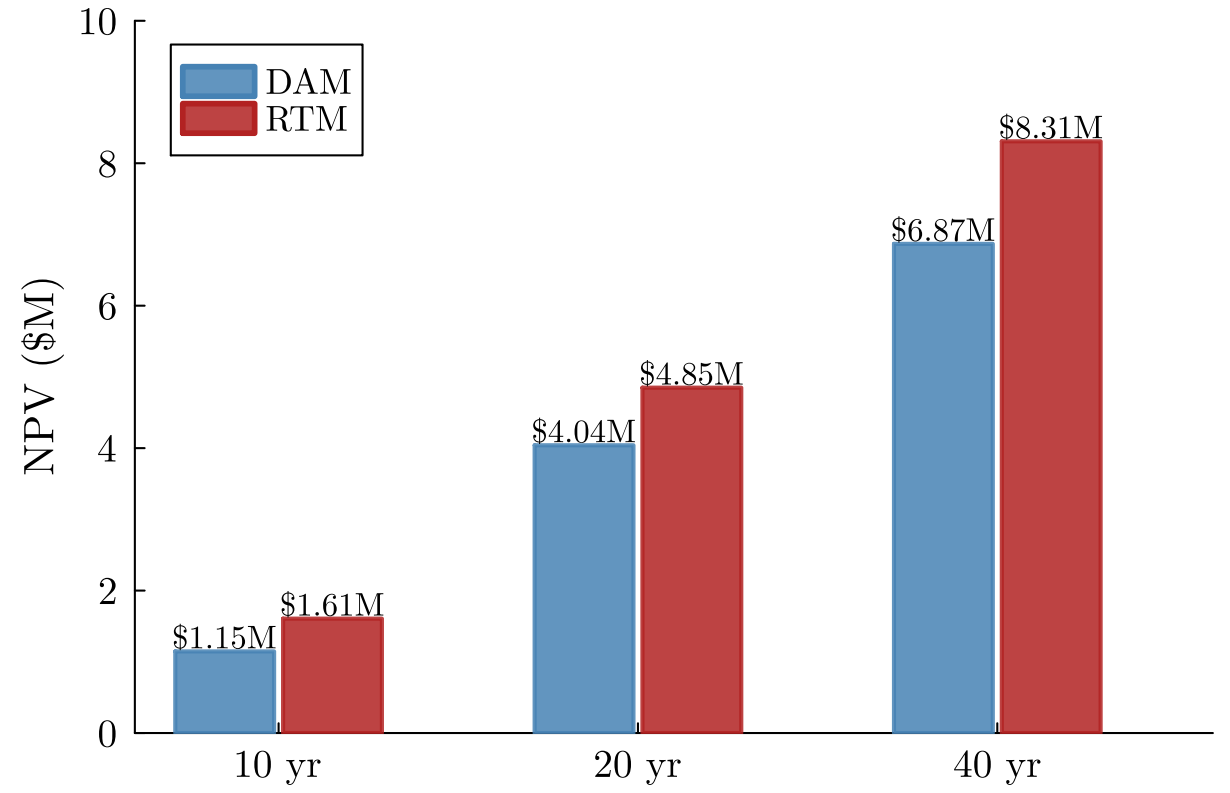

**Figure 3.** Net present value comparison between the RTM and the DAM for a 2.2 MW AWE across varying time horizons from Table 2. Values derived from final UB solution of $c = 24$ aggregate-informed Benders.

This algorithm development now enables us to assess the performance of electrolysis devices participating in either the DAM or the RTM. Figure 3 shows the NPV of market participation in those markets across the 10, 20, 40-year time horizons using the solutions from our algorithm. We find that participation in the RTM exceeds the DAM with $460k, $800k, and $1.44M NPV differ-

ences for each time horizon, respectively. We expect this due to the increased volatility in the intraday markets which enables taking advantage of market arbitrage more often, coinciding with previous studies [2]. In-depth analysis of these results and further market participation comparisons is left for future work.

## CONCLUSIONS

In this work, we introduced an aggregate-informed Benders decomposition algorithm for solving a multi-scale electrolysis system market participation optimization model. We empirically showed the superior capabilities of the proposed algorithm for solving large-scale market participation models. We expanded on previous work to compare economic performance between the DAM and RTM. We demonstrated that it remains more economically beneficial to participate in the RTM due to the volatility in intraday markets. For future work, we will explore the inherent trade-off of granularity in the aggregate and tractability, along with exploring the impact of clustering procedures. Additionally, we will explore the extent to which greedy solutions approach the monolithic solution in the upper bounding procedure.

## ACKNOWLEDGEMENTS

We acknowledge financial support from the U.S. National Science Foundation (CMMI-2328160).

## DELCARATION OF USE OF AI

Generative AI was used in programming the clustering algorithm and in producing the figures. The authors reviewed and edited the content as needed and take full responsibility for the content of the published article.

## AUTHOR IDENTIFIERS

Author ORCIDs:
Jennings K.X. 0009-0001-4007-1829
Zavala V.M. 0000-0002-5744-7378
Avraamidou S. 0000-0002-9334-9951

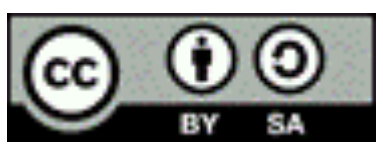